\documentclass[12pt,a4paper,intlimits,reqno]{amsart}
\usepackage[a4paper,top=4cm,bottom=3.5cm,left=3cm,right=3cm]{geometry}
\usepackage{amsmath, amssymb, mathrsfs, amsthm, enumerate}

\theoremstyle{plain}
\newtheorem{thm}{Theorem}

\theoremstyle{definition}
\newtheorem{exmp}{Example}
\newtheorem{rem}{Remark}

\DeclareMathOperator{\supp}{supp}

\title[Cyclic vectors of self-adjoint operators in Hilbert space]{Cyclic vectors of self-adjoint operators \\ in Hilbert space}
\author{Hidayat M. Huseynov}
\address{Department of Applied Mathematics, Baku State University, 23 Z.Khalilov str., AZ1148, Baku, Azerbaijan}
\address{Institute of Mathematics and Mechanics, National Academy of Sciences of Azerbaijan, 9 F.Agayev str., AZ1141, Baku, Azerbaijan}
\email{hmhuseynov@gmail.com}

\begin{document}
\maketitle

\section*{Introduction}

Let $H$ be a separable Hilbert space with inner product $(.,.)$ and norm $\| \cdot \| = \sqrt{(.,.)}$, and let $A$ be a bounded or unbounded self-adjoint operator in this space with a domain $D(A)$. A vector $f \in C^{\infty}(A) = \bigcap_{n=1}^{\infty} D(A^n)$ is said to be a cyclic vector for $A$ if the closure $\mathcal{L}$ of the span of the vectors $f, Af, A^2f, \dots$ coincides with the space $H$, i.e. the system $\{f, Af, A^2f, \dots\}$ is complete in the space $H$ (\cite{AkhiezerGlazman}, \cite{ReedSimon}).

The problem of finding conditions for a vector to be a cyclic vector for the given operator is a hard problem, but for some concrete operators even the criteria for a vector to be a cyclic vector are found (see survey paper \cite{Nikolskii}).

In the present paper we obtain a criterion and sufficient conditions for a vector to be a cyclic vector for a class of self-adjoint operators, more precisely for self-adjoint operators $A$ satisfying the following condition:

a) the spectrum of the operator $A$ consists of simple eigenvalues $\lambda_j$: $\lambda_j < \lambda_{j+1}$, $j = 0, \pm 1, \pm 2, \dots$.

Then the system of eigenvectors $\{e_j\}_{j=-\infty}^{\infty}$ of the operator $A$ ($Ae_j = \lambda_j e_j$) forms an orthonormal basis in the space $H$.

Note that for self-adjoint operators the simplicity of the spectrum and the existence of at least one cyclic vector are equivalent (see, e.g., \cite{AkhiezerGlazman}).

In order to formulate the main results of this paper we introduce the following notation:
\begin{equation*}
  P_{2n+1}(\lambda) = \prod_{i=-n}^{n} (\lambda-\lambda_i), \quad \dot{P}_{2n+1}(\lambda) = \frac{d}{d\lambda} P_{2n+1}(\lambda);
\end{equation*}
$E_{2n+1}$ is a unit matrix in Euclidean space $\mathbb{C}_{2n+1}$; \\
$K_{2n+1}$ is a square matrix of order $(2n+1)$ with elements
\begin{equation*}
  k_{ij}^{(n)} = \frac{1}{(f, e_i)\overline{(f, e_j)}} \cdot \frac{1}{\dot{P}_{2n+1}(\lambda_i) \dot{P}_{2n+1}(\lambda_j)} \sum_{|s|>n} \frac{P_{2n+1}^2(\lambda_s) |(f, e_s)|^2}{(\lambda_s-\lambda_i)(\lambda_s-\lambda_j)};
\end{equation*}
\begin{equation*}
  \overline{k_{ij}^{(n)}} = k_{ji}^{(n)}, \quad i,j = -n, -n+1, \dots, n;
\end{equation*}
$\langle .,. \rangle$ is the inner product in Euclidean space $\mathbb{C}_{2n+1}$; \\
$e_{2n+1}^{(k)}$ is the $(2n+1)$-dimensional column vector with components $\delta_{ik}$, $i = -n, -n+1, \dots, n$, $|k| < n$. Here $\delta_{ik}$ is the Kronecker delta.

\begin{thm}
  Let the self-adjoint operator $A$ satisfy the property a). In order that $f \in C^{\infty}(A)$ be a cyclic vector for the operator $A$ it is necessary and sufficient that for each integer $k$ the following conditions hold:
\begin{align*}
1^{\circ}. \; & (f, e_k) \ne 0; \\
2^{\circ}. \; & \lim_{n \to \infty} \langle (E_{2n+1} + K_{2n+1})^{-1} e_{2n+1}^{(k)}, e_{2n+1}^{(k)} \rangle = 1.
\end{align*}
\end{thm}

\begin{thm}
  Let the self-adjoint operator $A$ satisfy the property a), and for each integer $k$ let $f \in C^{\infty}(A)$ satisfy the conditions
\begin{align*}
1^{\circ}. \; & (f, e_k) \ne 0; \\
2^{\circ}. \; & \lim_{n \to \infty} \frac{1}{\dot{P}_{2n+1}^2(\lambda_k)} \cdot \sum_{|s|>n} \frac{P_{2n+1}^2(\lambda_s) |(f, e_s)|^2}{(\lambda_s - \lambda_k)^2} = 0.
\end{align*}
Then the vector $f$ is a cyclic vector for the operator $A$.
\end{thm}

Applying Theorem 2, we prove the following theorem.

\begin{thm}
  Let there exists $C>0$ such that for all integers $k$ the Fourier coefficients of the $2\pi$-periodic function $f(x)$ satisfy the conditions
\begin{equation*}
  0 < \left| \int_{-\pi}^{\pi} f(x) e^{-ikx} dx \right| \le C e^{-\delta |k|},
\end{equation*}
where $\delta > (6c_0^2+2) / 3c_0^3$, $c_0$ is the positive solution of the equation $c^2 = e^{1/c^2}$. Then the system of successive derivatives of the function $f(x)$, i.e. the system of functions
\begin{equation*}
  f(x), f'(x), f''(x), \dots
\end{equation*}
forms a complete system in the space $L_2(-\pi, \pi)$.
\end{thm}

(Note that $c_0 = 1.328\dots$, $(6c_0^2+2) / 3c_0^3 = 1.79\dots$)

\begin{exmp}
  For the Fourier coefficients of the function $f(x) = e^{a \cos x}$ $(a \ne 0)$ we have (see, e.g., \cite[Ch.~2, \S 10]{Olver})
\begin{equation*}
  \int_{-\pi}^{\pi} f(x) e^{-ikx} dx = 2 \int_{0}^{\pi} e^{a \cos x} \cos kx dx = 2 \pi \left( \frac{a}{2} \right)^k \sum_{m=0}^{+\infty} \frac{\left| a/2 \right|^{2m}}{m!(m+k)!}.
\end{equation*}
It follows that
\begin{equation*}
  0 < \left| \int_{-\pi}^{\pi} f(x) e^{-ikx} dx \right| \le 2\pi \left| \frac{a}{2} \right|^k e^{a^2/4} \frac{1}{k!}.
\end{equation*}
Therefore the conditions of Theorem~3 hold and hence the system of functions
\begin{equation*}
  e^{a \cos x}, \left( e^{a \cos x} \right)', \left( e^{a \cos x} \right)'', \dots
\end{equation*}
is complete in the space $L_2(-\pi, \pi)$.
\end{exmp}

Systems of derivatives of an analytic function were considered as complete systems in spaces of analytic functions in \cite{Gromov}, \cite{Ibragimov}, \cite{Kazmin}, \cite{Markushevich} etc.

\section*{Proof of Theorem 1}

Let $f \in C^{\infty}(A)$ and assume that condition $1^{\circ}$ of Theorem~1 holds. We denote by $L_{2n+1}(f)$ the span of the vectors $f, Af, \dots, A^{2n}f$. Then the distance from the vector $e_k$ to the subspace $L_{2n+1}(f)$ is expressed as (\cite[p.~20]{AkhiezerGlazman})
\begin{equation*}
  \rho(e_k, L_{2n+1}(f)) = \sqrt{\frac{\Gamma(e_k, f, Af, \dots, A^{2n}f)}{\Gamma(f, Af, \dots, A^{2n}f)}},
\end{equation*}
where $\Gamma(g_1, g_2, \dots, g_m)$ is the Gram determinant of the vectors $g_1, g_2, \dots, g_m$:
\begin{equation*}
  \Gamma(g_1, g_2, \dots, g_m) = \det
\begin{pmatrix}
  (g_1, g_1) & (g_1, g_2) & \dots  & (g_1, g_m) \\
  (g_2, g_1) & (g_2, g_2) & \dots  & (g_2, g_m) \\
  \vdots     & \vdots     & \ddots & \vdots \\
  (g_m, g_1) & (g_m, g_2) & \dots  & (g_m, g_m) \\
\end{pmatrix}.
\end{equation*}
It is easy to show that
\begin{equation} \label{eq:rho^2}
  \rho^2(e_k, L_{2n+1}(f)) = 1 - \langle A_{2n+1}^{-1} b_{2n+1}^{(k)}, b_{2n+1}^{(k)}\rangle,
\end{equation}
where $A_{2n+1}$ is the Gram matrix and $b_{2n+1}^{(k)}$ is the vector from $\mathbb{C}_{2n+1}$:
\begin{equation} \label{eq:A_2n+1}
  A_{2n+1} =
\begin{pmatrix}
  (f, f)       & (f, Af)       & \dots  & (f, A^{2n}f) \\
  (Af, f)      & (Af, Af)      & \dots  & (Af, A^{2n}f) \\
  \vdots       & \vdots        & \ddots & \vdots \\
  (A^{2n}f, f) & (A^{2n}f, Af) & \dots  & (A^{2n}f, A^{2n}f)
\end{pmatrix},\;\;
  b_{2n+1}^{(k)} =
\begin{pmatrix}
  1 \\
  \lambda_k \\
  \vdots \\
  \lambda_k^{2n}
\end{pmatrix} (f,e_k).
\end{equation}
Indeed, if we denote $\vec{0} = (0, 0, \dots, 0)^T \in C_{2n+1}$, then
\begin{multline*}
  \rho^2(e_k, L_{2n+1}(f)) = \det \left\{
  \begin{pmatrix}
    1       & \vec{0}^{*} \\
    \vec{0} & A_{2n+1}
  \end{pmatrix}^{-1}
  \begin{pmatrix}
    1              & b_{2n+1}^{(k)^{*}} \\
    b_{2n+1}^{(k)} & A_{2n+1}
  \end{pmatrix}
\right\} \\
= \det
  \begin{pmatrix}
    1                            & b_{2n+1}^{(k)^{*}} \\
    A_{2n+1}^{-1} b_{2n+1}^{(k)} & E_{2n+1}
  \end{pmatrix}
= 1 - \langle A_{2n+1}^{-1} b_{2n+1}^{(k)}, b_{2n+1}^{(k)}\rangle.
\end{multline*}
Using the eigenvector expansion for the self-adjoint operator $A^{\ell}$,
\begin{equation*}
  A^{\ell} f = \sum_{i=-\infty}^{\infty} \lambda_i^{\ell} (f, e_i) e_i, \quad \ell = 0, 1, 2, \dots,
\end{equation*}
it is easy to check that
\begin{equation} \label{eq:A_2n+1_by_B}
  A_{2n+1} = \sum_{j=-\infty}^{+\infty} B_{2n+1,j} B_{2n+1,j}^{*},
\end{equation}
where $B^{*}$ denotes the adjoint matrix of $B$ and the matrices $B_{s,j}$ $(s=2n+1)$ have the form
\begin{equation} \label{eq:B_sj}
  B_{s,j} =
\begin{pmatrix}
  (f, e_{sj-n})                     & (f, e_{sj-n+1})                       & \dots  & (f, e_{sj+n}) \\
  \lambda_{sj-n}(f, e_{sj-n})       & \lambda_{sj-n+1}(f, e_{sj-n+1})       & \dots  & \lambda_{sj+n}(f, e_{sj+n}) \\
  \vdots                            & \vdots                                & \ddots & \vdots \\
  \lambda_{sj-n}^{2n} (f, e_{sj-n}) & \lambda_{sj-n+1}^{2n} (f, e_{sj-n+1}) & \dots  & \lambda_{sj+n}^{2n} (f, e_{sj+n})
\end{pmatrix}.
\end{equation}
In the sequel $k$ is a fixed integer and $n$ is a natural number such that $|k|<n$. We put $c_{2n+1}^{(k)} = A_{2n+1}^{-1} b_{2n+1}^{(k)}$. Then $A_{2n+1} c_{2n+1}^{(k)} = b_{2n+1}^{(k)}$ and according to~\eqref{eq:A_2n+1_by_B}
\begin{equation*}
  \sum_{j=-\infty}^{+\infty} B_{2n+1,j} B_{2n+1,j}^{*} c_{2n+1}^{(k)} = b_{2n+1}^{(k)},
\end{equation*}
or
\begin{equation*}
  B_{2n+1,0} B_{2n+1,0}^{*} c_{2n+1}^{(k)} + \sum_{|j|>0} B_{2n+1,j} B_{2n+1,j}^{*} c_{2n+1}^{(k)} = b_{2n+1}^{(k)}.
\end{equation*}
If we put $\widehat{c}_{2n+1}^{(k)} = B_{2n+1,0}^{*} c_{2n+1}^{(k)}$ we obtain
\begin{equation} \label{eq:c_2n+1^k}
  \widehat{c}_{2n+1}^{(k)} + K_{2n+1} \widehat{c}_{2n+1}^{(k)} = B_{2n+1,0}^{-1} b_{2n+1}^{(k)},
\end{equation}
where the self-adjoint matrix $K_{2n+1}$ has the form
\begin{equation} \label{eq:K_2n+1}
  K_{2n+1} = \sum_{|j|>0} B_{2n+1,0}^{-1} B_{2n+1,j} B_{2n+1,j}^{*} B_{2n+1,0}^{*-1}.
\end{equation}
According to~\eqref{eq:A_2n+1} and~\eqref{eq:B_sj}
\begin{equation*}
  B_{2n+1,0}^{-1} b_{2n+1}^{(k)} = e_{2n+1}^{(k)},
\end{equation*}
where $e_{2n+1}^{(k)}$ is the column vector from $C_{2n+1}$ with components $\delta_{ik}$, $i = -n, -n+1, \dots, n$. Therefore the equation~\eqref{eq:c_2n+1^k} can be written in the form
\begin{equation*}
  \widehat{c}_{2n+1}^{(k)} + K_{2n+1} \widehat{c}_{2n+1}^{(k)} = e_{2n+1}^{(k)}.
\end{equation*}
From here, we have
\begin{equation*}
  \widehat{c}_{2n+1}^{(k)} = \left( E_{2n+1} + K_{2n+1} \right)^{-1} e_{2n+1}^{(k)}.
\end{equation*}
Now we can express the right-hand side of the formula~\eqref{eq:rho^2} in terms of the matrix $K_{2n+1}$ and the vector $e_{2n+1}^{(k)}$:
\begin{multline*}
  \rho^2(e_k, L_{2n+1}(f)) = 1 - \langle A_{2n+1}^{-1} b_{2n+1}^{(k)}, b_{2n+1}^{(k)}\rangle = 1 - \langle c_{2n+1}^{(k)}, b_{2n+1}^{(k)}\rangle \\
  = 1 - \langle B_{2n+1,0}^{*-1} \widehat{c}_{2n+1}^{(k)}, b_{2n+1}^{(k)}\rangle = 1 - \langle \widehat{c}_{2n+1}^{(k)}, B_{2n+1,0}^{-1} b_{2n+1}^{(k)}\rangle \\
  = 1 - \langle \left( E_{2n+1} + K_{2n+1} \right)^{-1} e_{2n+1}^{(k)}, e_{2n+1}^{(k)}\rangle.
\end{multline*}
So we obtain the following main formula which will play an important role in the sequel:
\begin{equation} \label{eq:rho^2_in_terms_of_K_and_e}
  \rho^2(e_k, L_{2n+1}(f)) = 1 - \langle \left( E_{2n+1} + K_{2n+1} \right)^{-1} e_{2n+1}^{(k)}, e_{2n+1}^{(k)}\rangle.
\end{equation}
We express the elements of the matrix $K_{2n+1}$ in terms of the eigenvalues $\lambda_j$ and the Fourier coefficients $(f, e_j)$ of the element $f$. First we consider the matrix $B_{2n+1,0}^{-1} B_{2n+1,j}$. Using the form~\eqref{eq:B_sj} for the matrix $B_{s,j}$ we can write
\begin{multline*}
  B_{2n+1,0}^{-1} =
\begin{pmatrix}
  (f, e_{-n})                   & (f, e_{-n+1})                    & \dots  & (f, e_{n}) \\
  \lambda_{-n}(f, e_{-n})       & \lambda_{-n+1}(f, e_{-n+1})      & \dots  & \lambda_{n}(f, e_{n}) \\
  \vdots                        & \vdots                           & \ddots & \vdots \\
  \lambda_{-n}^{2n} (f, e_{-n}) & \lambda_{-n+1}^{2n} (f, e_{-n+1}) & \dots  & \lambda_{n}^{2n} (f, e_{n})
\end{pmatrix}^{-1} \\
  = \frac{1}{\det B_{2n+1,0}}
\begin{pmatrix}
  \widehat{B}_{-n,-n}   & \widehat{B}_{-n+1,-n}   & \dots  & \widehat{B}_{n,-n} \\
  \widehat{B}_{-n,-n+1} & \widehat{B}_{-n+1,-n+1} & \dots  & \widehat{B}_{n,-n+1} \\
  \vdots                & \vdots                  & \ddots & \vdots \\
  \widehat{B}_{-n,n}    & \widehat{B}_{-n+1,n}    & \dots  & \widehat{B}_{n,n}
\end{pmatrix},
\end{multline*}
where $\det B_{2n+1,0} = \prod_{j=-n}^{n} (f, e_j) W(\lambda_{-n}, \dots, \lambda_{n})$, $W(\lambda_{-n}, \dots, \lambda_{n})$ is the Vandermonde determinant
\begin{equation} \label{eq:W}
  W(\lambda_{-n}, \dots, \lambda_{n}) =
\begin{pmatrix}
  1                 & 1                   & \dots  & 1 \\
  \lambda_{-n}      & \lambda_{-n+1}      & \dots  & \lambda_{n} \\
  \vdots            & \vdots              & \ddots & \vdots \\
  \lambda_{-n}^{2n} & \lambda_{-n+1}^{2n} & \dots  & \lambda_{n}^{2n}
\end{pmatrix},
\end{equation}
and $\widehat{B}_{ij}$ are the algebraic complements of the elements of the matrix $B_{2n+1,0}$. Therefore
\begin{multline*}
  B_{2n+1,0}^{-1} B_{2n+1,j} = \frac{1}{\prod_{j=-n}^{n} (f, e_j) W(\lambda_{-n}, \dots, \lambda_{n})} \\
\times \begin{pmatrix}
  \sum\limits_{s=-n}^{n} \widehat{B}_{s,-n} \lambda_{(2n+1)j-n}^{s+n} (f, e_{(2n+1)j-n}) & \dots  & \sum\limits_{s=-n}^{n} \widehat{B}_{s,-n} \lambda_{(2n+1)j+n}^{s+n} (f, e_{(2n+1)j+n}) \\
  \vdots                                                                          & \ddots & \vdots \\
  \sum\limits_{s=-n}^{n} \widehat{B}_{s,n} \lambda_{(2n+1)j-n}^{s+n} (f, e_{(2n+1)j-n})  & \dots  & \sum\limits_{s=-n}^{n} \widehat{B}_{s,n} \lambda_{(2n+1)j+n}^{s+n} (f, e_{(2n+1)j+n})
\end{pmatrix}.
\end{multline*}
Since
\begin{equation*}
  \sum_{s=-n}^{n} \widehat{B}_{s,i} \mu^{s+n} (f, e) = \prod_{s=-n \atop s \ne i}^{n} (f, e_s) (f, e) W(\lambda_{-n}, \dots, \mu, \dots, \lambda_{n}),
\end{equation*}
where $W(\lambda_{-n}, \dots, \mu, \dots, \lambda_{n})$ denotes the Vandermonde determinant which is obtained from~\eqref{eq:W} by replacing its $i$-th column by the vector $(1, \mu, \mu^2, \dots, \mu^{2n})^T$, we have
\begin{multline*}
  \frac{\sum_{s=-n}^{n} \widehat{B}_{s,i} \mu^{s+n} (f, e)}{\det B_{2n+1,0}} = \frac{\prod_{s=-n \atop s \ne i}^{n} (f, e_s) (f, e) W(\lambda_{-n}, \dots, \mu, \dots, \lambda_{n})}{\prod_{j=-n}^{n} (f, e_j) W(\lambda_{-n}, \dots, \lambda_{n})} \\
  = \frac{(f, e)}{(f, e_i)} \frac{(\mu - \lambda_{-n}) \dots (\mu - \lambda_{i-1}) (\lambda_{i+1} - \mu) \dots (\lambda_{n} - \mu)}{(\lambda_i - \lambda_{-n}) \dots (\lambda_i - \lambda_{i-1}) (\lambda_{i+1} - \lambda_i) \dots (\lambda_{n} - \lambda_i)} \\
  = \frac{(f, e)}{(f, e_i)} \frac{P_{2n+1}(\mu)}{\dot{P}_{2n+1}(\lambda_i)} \frac{1}{\mu - \lambda_i}, \qquad \qquad i = -n, -n+1, \dots, n,
\end{multline*}
where $P_{2n+1}(\mu) = \prod_{i=-n}^{n} (\mu-\lambda_i)$ is the polynomial of degree $(2n+1)$ and $\dot{P}_{2n+1}(\mu) = \frac{d}{d\mu} P_{2n+1}(\mu)$. Therefore
\begin{multline*}
  B_{2n+1,0}^{-1} B_{2n+1,j} \\
= \begin{pmatrix}
  \frac{(f, e_{(2n+1)j-n}) P_{2n+1}(\lambda_{(2n+1)j-n})}{(f, e_{-n}) \dot{P}_{2n+1}(\lambda_{-n}) (\lambda_{(2n+1)j-n} - \lambda_{-n})} & \dots  & \frac{(f, e_{(2n+1)j+n}) P_{2n+1}(\lambda_{(2n+1)j+n})}{(f, e_{-n}) \dot{P}_{2n+1}(\lambda_{-n}) (\lambda_{(2n+1)j+n} - \lambda_{-n})} \\
  \vdots                                                                          & \ddots & \vdots \\
  \frac{(f, e_{(2n+1)j-n}) P_{2n+1}(\lambda_{(2n+1)j-n})}{(f, e_{n}) \dot{P}_{2n+1}(\lambda_{n}) (\lambda_{(2n+1)j-n} - \lambda_{n})} & \dots  & \frac{(f, e_{(2n+1)j+n}) P_{2n+1}(\lambda_{(2n+1)j+n})}{(f, e_{n}) \dot{P}_{2n+1}(\lambda_{n}) (\lambda_{(2n+1)j+n} - \lambda_{n})}
\end{pmatrix}.
\end{multline*}
Then from~\eqref{eq:K_2n+1} we have the following expression for the elements $k_{ij}^{(n)}$ of the matrix $K_{2n+1}$:
\begin{equation} \label{eq:k_ij}
  k_{ij}^{(n)} = \frac{1}{(f, e_i)\overline{(f, e_j)}} \cdot \frac{1}{\dot{P}_{2n+1}(\lambda_i) \dot{P}_{2n+1}(\lambda_j)} \sum_{|s|>n} \frac{P_{2n+1}^2(\lambda_s) |(f, e_s)|^2}{(\lambda_s-\lambda_i)(\lambda_s-\lambda_j)},
\end{equation}
\begin{equation*}
  \overline{k_{ij}^{(n)}} = k_{ji}^{(n)}, \quad i,j = -n, -n+1, \dots, n.
\end{equation*}

Now, if all the conditions of Theorem 1 hold, then it follows from~\eqref{eq:rho^2_in_terms_of_K_and_e} that $\lim_{n \to \infty} \rho(e_k, L_{2n+1}(f)) = 0$ for any fixed integer $k$, i.e. all eigenvectors $e_k$ of the operator $A$ belong to the closure $\mathcal{L}$ of the span of the vectors $f, Af, A^2f, \dots$.

Indeed, according to the definition of the limit, for any $\varepsilon > 0$, we can find $n_{\varepsilon}$ such that
\begin{equation*}
  \alpha = \inf_{c_i} \| e_k - \sum_{i=0}^{n_{\varepsilon}} c_i A^i f \| < \varepsilon.
\end{equation*}
By the definition of the infimum, we also have that there are $c_i^{\varepsilon}$ such that
\begin{equation*}
  \| e_k - \sum_{i=0}^{n_{\varepsilon}} c_i^{\varepsilon} A^i f \| < \alpha + \varepsilon < 2 \varepsilon.
\end{equation*}
Since $\{e_k\}$ forms an orthonormal basis in the Hilbert space $H$ and $e_k \in \mathcal{L}$, we have $\mathcal{L} = H$.

The necessity of conditions $1^{\circ}$ and $2^{\circ}$ for the equality $\mathcal{L} = H$ is obvious.
\qed

\section*{Proof of Theorem 2}

It is sufficient to prove the inequality ($|k|<n$)
\begin{equation} \label{eq:rho^2_inequality}
  \rho^2(e_k, L_{2n+1}(f)) \le \frac{1}{|(f, e_k)|^2} \cdot \frac{1}{\dot{P}_{2n+1}^2(\lambda_k)} \sum_{|s|>n} \frac{P_{2n+1}^2(\lambda_s) |(f, e_s)|^2}{(\lambda_s-\lambda_k)^2}.
\end{equation}

We use the following easy fact (see, e.g., \cite{MarcusMinc}): if $M$ is a positive definite Hermitian square matrix of order $m \times m$, $\langle .,. \rangle_m$ is the inner product in Euclidean space $\mathbb{C}_{m}$ and $e \in \mathbb{C}_{m}$, $\|e\|_m = 1$, then the following inequality holds:
\begin{equation} \label{eq:MarcusMinc}
  \langle M^{-1} e, e \rangle_m \cdot \langle M e, e \rangle_m \ge 1.
\end{equation}

If we put $M = E_{2n+1} + K_{2n+1}$, $m = 2n+1$, $e = e_{2n+1}^{(k)}$ from~\eqref{eq:MarcusMinc} we obtain
\begin{equation*}
  \langle (E_{2n+1} + K_{2n+1})^{-1} e_{2n+1}^{(k)}, e_{2n+1}^{(k)} \rangle \ge \frac{1}{1 + \langle K_{2n+1} e_{2n+1}^{(k)}, e_{2n+1}^{(k)} \rangle}.
\end{equation*}
Then from~\eqref{eq:rho^2_in_terms_of_K_and_e} we have
\begin{equation} \label{eq:rho^2_inequality_2}
  \rho^2(e_k, L_{2n+1}(f)) \le \frac{\langle K_{2n+1} e_{2n+1}^{(k)}, e_{2n+1}^{(k)} \rangle}{1 + \langle K_{2n+1} e_{2n+1}^{(k)}, e_{2n+1}^{(k)} \rangle}.
\end{equation}
On the other hand, according to~\eqref{eq:k_ij}
\begin{equation} \label{eq:K_e_e}
  \langle K_{2n+1} e_{2n+1}^{(k)}, e_{2n+1}^{(k)} \rangle = k_{kk}^{(n)} = \frac{1}{|(f, e_k)|^2} \cdot \frac{1}{\dot{P}_{2n+1}^2(\lambda_k)} \sum_{|s|>n} \frac{P_{2n+1}^2(\lambda_s) |(f, e_s)|^2}{(\lambda_s-\lambda_k)^2}.
\end{equation}
Relations~\eqref{eq:rho^2_inequality_2} and~\eqref{eq:K_e_e} imply inequality~\eqref{eq:rho^2_inequality}. \qed

\begin{rem}
  In the proof of~\eqref{eq:rho^2_inequality}, we assumed that all Fourier coefficients $(f, e_s)$ of the vector $f$ differ from zero. Now consider the case when some of the Fourier coefficients are equal to zero. Denote by $(f, e_{m_s})$, $s = 0, \pm 1, \pm 2, \dots$ all nonzero Fourier coefficients. Then eigenvector expansion formula for the vector $f$ has the form
\begin{equation*}
  A^{\ell} f = \sum_{i=-\infty}^{\infty} \widetilde{\lambda}_i^{\ell} (f, \widetilde{e}_i) \widetilde{e}_i, \quad \ell = 0, 1, 2, \dots, \quad \textrm{ where } \widetilde{\lambda}_s = \lambda_{m_s}, \; \widetilde{e}_s = e_{m_s}.
\end{equation*}
Repeating the arguments used in the proofs of Theorems 1 and 2 we obtain the following analogue of inequality~\eqref{eq:rho^2_inequality} ($|k| < n$):
\begin{equation*}
  \rho^2(\widetilde{e}_k, L_{2n+1}(f)) \le \frac{1}{|(f, \widetilde{e}_k)|^2} \cdot \frac{1}{\dot{\widetilde{P}}_{2n+1}^2(\widetilde{\lambda}_k)} \sum_{|s|>n} \frac{\widetilde{P}_{2n+1}^2(\widetilde{\lambda}_s) |(f, \widetilde{e}_s)|^2}{(\widetilde{\lambda}_s-\widetilde{\lambda}_k)^2},
\end{equation*}
where
\begin{equation*}
  \widetilde{P}_{2n+1}(\lambda) = \prod_{i=-n}^{n} (\lambda-\widetilde{\lambda}_i).
\end{equation*}
\end{rem}

\begin{rem}
  In the proof of Theorem 2 we actually obtained the following estimate:
\begin{equation*}
  \rho^2(e_k, L_{2n+1}(f)) \le \frac{k_{kk}^{(n)}}{1 + k_{kk}^{(n)}}.
\end{equation*}

It can be shown that a more accurate estimate holds:
\begin{equation*}
  \rho^2(e_k, L_{2n+1}(f)) \le \frac{k_{kk}^{(n)}}{1 + k_{kk}^{(n)}} - \frac{\left( \sum\limits_{i=1 \atop i \ne k}^{n} |k_{ik}^{(n)}|^2 \right)^2}{\sum\limits_{i=1 \atop i \ne k}^{n} |k_{ik}^{(n)}|^2 + \sum\limits_{j=1 \atop j \ne k}^{n} \sum\limits_{i=1 \atop i \ne k}^{n} k_{kj}^{(n)} k_{ji}^{(n)} k_{ik}^{(n)}} \cdot \frac{1}{\left( 1 + k_{kk}^{(n)} \right)^2}.
\end{equation*}
\end{rem}

\section*{Proof of Theorem 3}

Consider the operator $A$ generated by the differential expression $i \, d/dx$ in the space $L_2(-\pi,\pi)$:
\begin{align*}
  Ay   & = i\frac{dy}{dx}, \\
  D(A) & = \{ y \in L_2(-\pi,\pi) | y \in \mathscr{AC}(-\pi,\pi), y(-\pi) = y(\pi), y' \in L_2(-\pi,\pi) \}.
\end{align*}

It is known (see, e.g., \cite{AkhiezerGlazman}) that $A$ is a self-adjoint operator in the Hilbert space $H = L_2(-\pi,\pi)$, $\lambda_k = k$ $(k = 0, \pm 1, \pm 2, \dots)$ are the eigenvalues and $e_k = \frac{1}{\sqrt{2\pi}} e^{-ikx}$ are the eigenfunctions of the operator $A$. So the operator $A$ satisfies the property~a). Let the function $f(x)$ satisfy the conditions of Theorem 3. Let us show that all the conditions of Theorem 2 are satisfied. From the estimate for the Fourier coefficients we obtain that $f(x)$ is an analytic function on the segment $[-\pi,\pi]$ (see, e.g., \cite[p.~90]{Bari}) and since this function is $2\pi$-periodic we have $f^{(m)}(-\pi) = f^{(m)}(\pi)$, $m = 0, 1, 2, \dots$, i.e. $f \in C^{\infty}(A)$. Now it is sufficient to prove the fulfilment of condition~$2^{\circ}$ of Theorem~2.

Since $P_{2n+1}(\lambda) = \prod_{i=-n}^{n} (\lambda-\lambda_i) = \prod_{i=-n}^{n} (\lambda-i) = \lambda (\lambda^2-1^2) \ldots (\lambda^2-n^2),$ then using the inequality $1 - x^2 < e^{-x^2}$, $0 < x < 1$ for $s > n$ we obtain
\begin{multline*}
  P_{2n+1}(\lambda_s) = s (s^2-1^2) \ldots (s^2-n^2) = s^{2n+1} \left( 1 - \left( \frac{1}{s} \right)^2 \right) \ldots \left( 1 - \left( \frac{n}{s} \right)^2 \right) \\
  < s^{2n+1} e^{- (1^2 + \ldots + n^2) / s^2} = s^{2n+1} e^{- n(n+1)(2n+1) / 6s^2} \le s^{2n+1} e^{- n^3 / 3s^2}.
\end{multline*}
Finally, taking into account that $|(f,e_s)| \le C e^{-\delta |s|}$ and
\begin{equation*}
  \dot{P}_{2n+1}(\lambda_k) = \prod_{i=-n \atop i \ne k}^{n} (k-i) = (-1)^{n-k} (n+k)! (n-k)!, \quad |k|<n
\end{equation*}
and putting $\delta = \sigma + \delta_1$, where $\sigma = (6c_0^2+2) / 3c_0^3$, $\delta_1 > 0$ we have
\begin{align*}
  \frac{1}{\dot{P}_{2n+1}^2(\lambda_k)} & \sum_{|s|>n} \frac{P_{2n+1}^2(\lambda_s) |(f, e_s)|^2}{(\lambda_s - \lambda_k)^2} \\
  & \le \frac{2 C^2}{\left( (n+k)! (n-k)! \right)^2} \sum_{s=n+1}^{+\infty} \frac{s^{4n+2} e^{- 2n^3 / 3s^2} e^{-2 \delta s}}{(s-k)^2} \\
  & \le \frac{2 C^2}{\left( (n+k)! (n-k)! \right)^2} \max_{n < s < \infty} \left( s^{4n} e^{- 2n^3 / 3s^2} e^{-2 \sigma s} \right) \sum_{s=n+1}^{+\infty} \frac{s^2 e^{-2 \delta_1 s}}{(s-k)^2} \\
  & = \frac{2 C^2}{\left( (n+k)! (n-k)! \right)^2} \left( s^{4n} e^{- 2n^3 / 3s^2} e^{-2 \sigma s} \right)_{s=c_0 n} \sum_{s=n+1}^{+\infty} \frac{s^2 e^{-2 \delta_1 s}}{(s-k)^2} \\
  & = \frac{2 C^2}{\left( (n+k)! (n-k)! \right)^2} (c_0 n)^{4n} e^{- 2n / 3c_0^2} e^{-2 \sigma c_0 n} \sum_{s=n+1}^{+\infty} \frac{s^2 e^{-2 \delta_1 s}}{(s-k)^2} \\
  & \le \frac{2 C^2 n^{4n} e^{-4n}}{\left( (n+k)! (n-k)! \right)^2} \frac{1}{(n+1-k)^2} e^{-2} \frac{e^{-\delta_1 (n+1)}}{1 - e^{-\delta_1}} \left( \frac{2}{\delta_1} \right)^2.
\end{align*}
From here it follows that for each fixed $k$
\begin{equation*}
  \lim_{n \to \infty} \frac{1}{\dot{P}_{2n+1}^2(\lambda_k)} \cdot \sum_{|s|>n} \frac{P_{2n+1}^2(\lambda_s) |(f, e_s)|^2}{(\lambda_s - \lambda_k)^2} = 0,
\end{equation*}
i.e. condition $2^{\circ}$ of Theorem~2 holds.
\qed

\begin{rem}
  Theorem 3 shows that the function $f(x)$ satisfying the conditions of this theorem is a cyclic vector for the self-adjoint operator $A$ generated in the space $L_2(-\pi,\pi)$ by the boundary value problem
\begin{align*}
  iy'     & = \lambda y, \\
  y(-\pi) & = y(\pi).
\end{align*}

Let us show that an infinitely differentiable finite function $\varphi(x)$ with support $\supp \varphi(x) \subset [-\pi,\pi]$ is not a cyclic vector for the operator $A$. Assume the contrary, i.e., assume that such a function is a cyclic vector for the operator $A$.

Consider the Fourier transform of the function $\varphi(x)$:
\begin{equation*}
  \Phi(\lambda) = \int_{-\pi}^{\pi} \varphi(x) e^{-i \lambda x} dx.
\end{equation*}
It is known (see, e.g., \cite[p.~22]{Levin}) that the function $\Phi(\lambda)$ has an infinite number of zeros. Let $\lambda_0$ be a zero of this function:
\begin{equation*}
  \int_{-\pi}^{\pi} \varphi(x) e^{-i \lambda_0 x} dx = 0.
\end{equation*}
Integrating by parts we obtain
\begin{equation*}
  \int_{-\pi}^{\pi} \varphi^{(j)}(x) e^{-i \lambda_0 x} dx = 0, \qquad j = 0, 1, 2, \ldots.
\end{equation*}
But this is impossible since by assumption the system of functions $\{ \varphi^{(j)}(x) \}_{j=0}^{+\infty}$ is complete in the space $L_2(-\pi,\pi)$.
\end{rem}

\end{document}